\documentclass[12pt, a4paper, reqno]{amsart}

\usepackage{amssymb}
\usepackage[all]{xy}
\usepackage[
  margin=1.5cm,
  includefoot,
  footskip=30pt,
  ]{geometry}
\usepackage[colorlinks=true,citecolor=cyan,linkcolor=magenta, urlcolor=blue, filecolor=green, backref=page]{hyperref}

\newcommand{\cO}{\mathcal{O}}

\newcommand{\bF}{\mathbb{F}}

\newcommand{\bQ}{\mathbb{Q}}

\newcommand{\bZ}{\mathbb{Z}}

\newcommand{\fp}{\mathfrak{p}}

\DeclareMathOperator{\Div}{Div}

\DeclareMathOperator{\Gal}{Gal}

\DeclareMathOperator{\Norm}{Norm}

\DeclareMathOperator{\Sel}{Sel}

\DeclareMathOperator{\im}{im}

\DeclareMathOperator{\rk}{rank}

\theoremstyle{plain}% default style
\newtheorem{theorem}{Theorem}[section]

\newtheorem{lemma}[theorem]{Lemma}
\newtheorem{proposition}[theorem]{Proposition}
\newtheorem{corollary}[theorem]{Corollary}

\theoremstyle{definition} % definition style

\theoremstyle{remark} % remark style
\newtheorem{remark}[theorem]{Remark}
\newtheorem{example}[theorem]{Example}

\DeclareFontFamily{U}{wncy}{}
\DeclareFontShape{U}{wncy}{m}{n}{<->wncyr10}{}
\DeclareSymbolFont{mcy}{U}{wncy}{m}{n}
\DeclareMathSymbol{\Sha}{\mathord}{mcy}{"58}

\newcommand{\lp}{\left(}
\newcommand{\rp}{\right)}

\newcommand{\lara}[1]{\langle #1 \rangle}
\newcommand{\lbrb}[1]{\lp #1 \rp}
\newcommand{\lcrc}[1]{\left\{ #1 \right\}}

\newcommand{\quadsym}[2]{\lbrb{\frac{#1}{#2}}}

\title{On the Jacobian of hyperelliptic curves $y^2 = x^5 + m^2$}
\author{Keunyoung Jeong}
\address{Department of Mathematical Sciences, Ulsan National Institute of Science and Technology, UNIST-gil 50, Ulsan 44919, Korea}
\email{kyjeongg@gmail.com}

\author{Junyeong Park}
\address{Department of Mathematics, POSTECH (Pohang University of Science and Technology), San 31, Hyoja-Dong, Nam-Gu, Pohang-si, Gyeongsangbuk-do 790-784, South Korea}
\email{junyeongp@gmail.com}

\author{Donggeon Yhee}
\address{Industrial and mathematical data analytics research center, Seoul National University, %Building 25 Room 128, 
Gwanak-Ro 1, Gwanak-Gu, Seoul, South Korea}
\email{dgyhee@gmail.com}

\begin{document}


\maketitle

\begin{abstract}
In this paper, we study the algebraic rank and the analytic rank of the Jacobian of hyperelliptic curves $y^2 = x^5 + m^2$ for integers $m$.
Namely, we first provide a condition on $m$ that gives a bound of the size of Selmer group and then we provide a condition on $m$ that makes $L$-functions non-vanishing.
As a consequence, we construct a Jacobian that satisfies the rank part of the Birch--Swinnerton-Dyer conjecture.
% On the Mordell--Weil group side, we give a condition on $m$ that makes the algebraic rank of Jacobian zero.
% On the $L$-value side, we give a sufficient condition on $m$ which guarantees the non-vanishing of $L$-values.
%As a consequence, we find infinitely many hyperelliptic curves whose Jacobian satisfies the rank part of the Birch--Swinnerton-Dyer conjecture. If we further assume the parity conjecture, then there are infinitely many Jacobians of hyperelliptic curve that satisfy Birch--Swinnerton-Dyer conjecture and have a non-trivial Shafarevich--Tate group.
\end{abstract}

\section{Introduction}

For each integer $A$, we define a hyperelliptic curve $C_A : y^2 = x^5 + A$ and its Jacobian $J_A$.
In \cite{Sto1, Sto2} Stoll studied the arithmetic of $C_A$ and in \cite{SY03} Stoll and Yang studied the $L$-values of $C_A$.
In this paper, we focus on the case of $A = m^2$ where $m$ is a square-free integer.
More precisely, we study the algebraic rank and the analytic rank of $J_{m^2}$.
We note that every hyperelliptic curve in our family does not satisfy the conditions \cite[(1.3)]{Sto1}, so this curve is not covered in \cite{Sto1}.

To get an algebraic rank, a standard method is to give a bound of the Selmer groups of the Jacobians. Using the result of Schaefer \cite{Sch} and the calculation of the root numbers \cite{Sto2}, we obtain the following.

\begin{theorem} \label{thm:main MW group}
There are infinitely many integers $m$ where $J = J_{m^2}$ satisfies
\begin{equation*}
    J(\bQ) \cong \bZ/5\bZ.
\end{equation*}
On the other hand, there are infinitely many $m$ such that
\begin{equation*}
J(\bQ) \cong \bZ/5\bZ \oplus \bZ
\end{equation*}
under the parity conjecture.
\end{theorem}
For simplicity, we mainly consider the case where $m$ is a prime.
However, our proof of this theorem can be applied to general $J_{m^2}$ for square-free $m$ such that all of the prime divisors $p$ of $m$ satisfy $p \not\equiv 1 \pmod{5}$, and there is at most one $p \equiv 4 \pmod{5}$ among them.
In this case, the primes of $K$ above $m$ satisfy a certain kind of orthogonality (i.e. there exist generators $\pi_\fp, \pi_{\fp'}$ such that $\pi_\fp$ is trivial in $K_{\fp'}^\times/K_{\fp'}^{\times 5}$ and vice versa).
This property makes the descent computation much easier as we will see in Proposition \ref{prop:descent not over 5}.
For the case where $m$ is not a prime, see Remark \ref{rem:} and Example \ref{exam:}.
As an example, we consider $m = 101$ a prime equivalent to 1 modulo 5 in Proposition \ref{prop:101}.

On the analytic side, there are results on the special $L$-value of the hyperelliptic curves $C_A$ like \cite{SY03, Mas}.
Such curves have complex multiplication, so there is a Hecke character $\eta_A$ satisfying
\begin{equation*}
    L(s, C_A) = L(s, J_A) = L(s, \eta_A).
\end{equation*}
Based on the work \cite{Yan97, Yan98, Yan99} on the non-vanishings of $L$-functions of Hecke characters and \cite{Sto1, Sto2} on hyperelliptic curves $C_{A}$, Stoll and Yang showed that
\begin{equation*}
    L(1, J_1) \neq 0
\end{equation*}
in \cite{SY03}.
In this paper, we extend this result for the curve $C_{A}$ with certain conditions on $A$, in Proposition \ref{prop:Lvalue} which gives an expression of $L(1, \eta_{A})$.
% In the expression of $L(1, \eta_A)$, the prime divisors of $A$ not equivalent to 1 modulo 5 make the Gauss? certain? sum (\ref{eqn:phiv1}), (\ref{eqn:phiv2}) which are not easy to see non-zero. However, in the case of prime divisors of $A$ equivalent to 1 modulo 5, the corresponding term would be non-zero constant.
As a consequence, we obtain
\begin{theorem} \label{thm:main2}
Let $J_{A}$ be a Jacobian of $C_{A}$ whose root number is $+1$.
If $A$ is a product of rational primes equivalent to 1 modulo 5 and $(A^4-1)$ is divided by 25, then $L(1, J_{A}) \neq 0$.
\end{theorem}

Note that the rational primes $p\equiv1\pmod5$ are exactly the ones split completely in $K$. In formula (\ref{eqn:Lvalue equation}), one can see from (\ref{eqn:Ivsplit}) that the factors involving primes $v$ of $F$ split in $K$ are non-zero. 
To see whether the factors involving primes of $F$ inert in $K$ vanish or not, one need to evaluate integral (\ref{eqn:defIv}), which seems to be complicated. However, when it comes to the descent on $C_{m^2}$, the situation seems complementary. More precisely, if $m$ only has prime factors which are not totally split, then the descent is manageable. However, if $m$ has prime factors which split completely in $K$, then the descent become more complicated to deal with. This explains why we cannot obtain an infinite family of Jacobians of the form $J_{m^2}$ satisfying the rank part of the Birch--Swinnerton-Dyer conjecture. Instead of this, we give an illustration for the case $p\equiv1\pmod5$:

% We can construct an example of abelian surface satisfying the rank part of the Birch--Swinnerton-Dyer conjecture when the algebraic rank zero condition and the analytic rank zero condition are satisfied simultaneously.
% However, our algebraic rank zero results and analytic rank zero results have no intersection.
% Hence we cannot make an infinite family of Jacobains satisfying the Birch-Swinnerton-Dyer conjecture. 
% Instead, we have
\begin{corollary} \label{cor:101}
A Jacobian $J_{101^2}$ satisfies the rank part of Birch--Swinnerton-Dyer conjecture.
\end{corollary}
We note that Corollary \ref{cor:101} may be deduced from 2-descent available in Magma and the numerical computation of $L$-values since the rank of $J_{101^2}$ is zero, but we want to emphasize that the analogous result for other primes $p \equiv 1 \pmod{5}$ may be deduced from our $(1- \zeta_5)$-descent with less computational complexity.

In Section \ref{sec:Preliminaries}, we list some facts on local fields and recall the computation of the root number of $J_{m^2}$.
Based on these results, we describe descent for Jacobians in Section \ref{sec:descent} and give a proof of Theorem \ref{thm:main MW group}.
After computing the special $L$-value in Section \ref{sec:L-value}, we will show Theorem \ref{thm:main2} and Corollary \ref{cor:101}.

\section{Preliminaries} \label{sec:Preliminaries}

\subsection{Local field computation}
We list some notations which will be used in Section \ref{sec:Preliminaries} and \ref{sec:descent}.
We fix a fifth root of unity $\zeta_5$ in $\overline{\bQ}$.
Let $K = \bQ(\zeta_5)$ and $F = \bQ(\sqrt{5})$.
We recall that a rational prime $p$ is inert, splits into two primes, splits completely in $K/\bQ$ if and only if $p \equiv 2$ or $3$, $p \equiv 4$, $p \equiv 1$ modulo 5, respectively.
In each case, we denote primes of $K$ above a rational prime $p$ by $p, w, v$ and its generator by $p, \pi_w, \pi_v$, respectively.
The unique prime above $5$ is denoted by $v_5$, but we also admit the notations $K_5$ and $\pi_5$ for $K_{v_5}$ and $\pi_{v_5}$.
We use a symbol $\fp$ to indicate a prime ideal of $K$ and $\pi$ to a prime element.
Also we use the notation $\zeta_n$ for a primitive $n$-th root unity in $K$ or any local fields, if it exists.

In this section, we compute the images of prime elements $\pi$ in $K_\fp^\times/K_\fp^{\times 5}$.
We first compute the group $K_\fp^{\times}/K_\fp^{\times n}$.
When $\fp = v_5$, we fix a generator $\pi_5$ by $(1 - \zeta_5)$. % for each embedding $K \to K_5$.
Since
\begin{equation*}
    K_5^\times \cong \pi_5^{\bZ} \times \mu_4 \times U^{(1)} \quad
    \textrm{and} \quad U^{(2)} \cong \bZ_5^4,
\end{equation*}
we have
\begin{equation} \label{eqn:local K5}
    K_5^\times/K_5^{\times 5} \cong \lara{\pi_5, 1+\pi_5, 1+ \pi_5^2, 1+ \pi_5^3, 1+ \pi_5^4, 1+ \pi_5^5}
\end{equation}
and every element in $U^{(6)}$ is a fifth-power.
We rename the generating elements by $ \lara{\alpha, \beta, \gamma, \delta, \epsilon, \eta}$.
For all other primes $\fp \neq v_5$, 5 is invertible in the ring of integers $\cO_{K, \fp}$.
 So we have
\begin{equation} \label{eqn: local field not 5}
    K_\fp^\times / K_\fp^{\times 5} \cong \lara{\pi_\fp, \zeta_{5^n}}
\end{equation}
where $\zeta_{5^n}$ generates the 5-part of the root of unities of $K_{\fp}^\times$.
We also rename the generating elements by $\lara{\alpha_\fp, \beta_\fp}$ and drop the subscript whenever the meaning is clear from the context.
We note that every element in $U^{(2)}$ is a fifth-power in this case.

We need $\pi_5$-expansions of some elements in $K_5$.
By expanding $\pi_5^4 = (1 - \zeta_5)^4$, we have
\begin{equation*}
    5 = 4\pi_5^4 + 3\pi_5^5 + 3 \pi_5^6 + 4\pi_5^7 + \pi_5^8  + 3 \pi_5^{9} + O(\pi_5^{11}).
\end{equation*}
We choose $\sqrt{5}$ and $\zeta_4$ in $K_5$ such that
\begin{align*}
    \sqrt{5}\equiv2\pi_5^2\pmod{\pi_5^3}
    \quad\textrm{and}\quad\zeta_4\equiv2\pmod{\pi_5}
\end{align*}
respectively. Then, one may verify that
\begin{align*}
    \sqrt{5} &= 2 \pi_5^2 + 2 \pi_5^3 + \pi_5^4 + O(\pi_5^7), \\
    \zeta_4 &= 2 + 4 \pi_5^4 + 3\pi_5^5 + O(\pi_5^6), \\
    \zeta_4^3 &= 3 + 2 \pi_5^4 + 4\pi_5^5 + O(\pi_5^6), \\
    -\left(\frac{1+\sqrt{5}}{2}\right) &=2+4\pi_5^2+4\pi_5^3+\pi_5^5+O(\pi_5^6),
\end{align*}
where the last one is a fundamental unit of $F$, which we will denote by $u_F$.
We note that $\{1, u_F \}$ is an integral basis of $\cO_F$, so we can choose a generator $\pi_w = a + b\sqrt{5}$ for $a, b \in \frac{1}{2}\bZ$, or $\pi_w = a + b u_F$ for $a, b \in \bZ$.

Now we can describe the images of the prime elements of $K$ which is not above a rational prime $p \equiv 1 \pmod{5}$ in $K_5^\times/K_5^{\times 5}$.

\begin{lemma} \label{lem:localover5}
(1) Let $n$ be a rational integer not divided by $5$. Then, the image of $n$ in $K_5^\times/K_5^{\times 5}$ is
\begin{equation*}
    \begin{array}{cl}
     1 &\textrm{if } n\equiv 1, 7, 18, 24 \pmod{25} \\
    \epsilon\eta^2 & \textrm{if } n\equiv 3,4, 21, 22 \pmod{25} \\
    \epsilon^2\eta^4 &\textrm{if } n\equiv 9, 12, 13, 16 \pmod{25} \\
    \epsilon^3\eta &\textrm{if } n\equiv 2, 11, 14, 23 \pmod{25} \\
    \epsilon^4\eta^3 &\textrm{if } n\equiv 6, 8, 17, 19\pmod{25}
    \end{array}
\end{equation*}
(2) For a prime $w$ above a rational prime $p \equiv 4 \pmod{5}$ and its generator $\pi_w =a + b\sqrt{5}$ with $a,b\in \frac{1}{2}\bZ$, the image of $\pi_w$ in $K_5^\times/K_5^{\times 5}$ is given by the following table.
\begin{equation*}
    \begin{array}{|c|c|c|c|c|c|} \hline
    a\pmod{5} & p\equiv 4 & p\equiv 9 & p\equiv 14 & p\equiv 19 & p\equiv 24\\
    \hline
    2 & \gamma^b\delta^b\epsilon^{b+3}\eta & \gamma^b\delta^b\epsilon^{b+1}\eta^2 & \gamma^b\delta^b\epsilon^{b+4}\eta^3 & \gamma^b\delta^b\epsilon^{b+2}\eta^4 & \gamma^b\delta^b\epsilon^{b}\\ \hline
%    4 & \gamma^{2b+4}\delta^{2b+4}\epsilon^{2b+2}\eta & \gamma^{2b+4}\delta^{2b+4}\epsilon^{2b}\eta^2 & \gamma^{2b+4}\delta^{2b+4}\epsilon^{2b+3}\eta^3 & \gamma^{2b+4}\delta^{2b+4}\epsilon^{2b+1}\eta^4 &\gamma^{2b+4}\delta^{2b+4}\epsilon^{2b+4}\\
    4 & \gamma^{3b}\delta^{3b}\epsilon^{3b+3}\eta & \gamma^{3b}\delta^{3b}\epsilon^{3b+1}\eta^2 &\gamma^{3b}\delta^{3b}\epsilon^{3b+4}\eta^3 & \gamma^{3b}\delta^{3b}\epsilon^{3b+2}\eta^4 &\gamma^{3b}\delta^{3b}\epsilon^{3b}\\ \hline
%    3 & \gamma^{4b}\delta^{4b}\epsilon^{4b+3}\eta & \gamma^{4b}\delta^{4b}\epsilon^{4b+1}\eta^2 & \gamma^{4b}\delta^{4b}\epsilon^{4b+4}\eta^3 & \gamma^{4b}\delta^{4b}\epsilon^{4b+2}\eta^4 & \gamma^{4b}\delta^{4b}\epsilon^{4b}\\
    3 & \gamma^{4b}\delta^{4b}\epsilon^{4b+3}\eta & \gamma^{4b}\delta^{4b}\epsilon^{4b+1}\eta^2 & \gamma^{4b}\delta^{4b}\epsilon^{4b+4}\eta^3 & \gamma^{4b}\delta^{4b}\epsilon^{4b+4}\eta^4 &\gamma^{4b}\delta^{4b}\epsilon^{4b}\\ \hline
%    1 & \gamma^{3b+4}\delta^{3b+4}\epsilon^{3b+2}\eta & \gamma^{3b+4}\delta^{3b+4}\epsilon^{3b}\eta^2 & \gamma^{3b+4}\delta^{3b+4}\epsilon^{3b+3}\eta^3 & \gamma^{3b+4}\delta^{3b+4}\epsilon^{3b+1}\eta^4 & \gamma^{3b+4}\delta^{3b+4}\epsilon^{3b+4}
    1 & \gamma^{2b}\delta^{2b}\epsilon^{2b+3}\eta & \gamma^{2b}\delta^{2b}\epsilon^{2b+1}\eta^2 & \gamma^{2b}\delta^{2b}\epsilon^{2b+4}\eta^3 & \gamma^{2b}\delta^{2b}\epsilon^{2b+1}\eta^4 & \gamma^{2b}\delta^{2b}\epsilon^{2b} 
    \\ \hline
    \end{array} 
\end{equation*}
\end{lemma}
\begin{proof}
For a generator $\sigma:\zeta_5\mapsto \zeta_5^2$ of $\Gal(K_5/\bQ_5)$, we have
\begin{align*}
    &\quad\sigma(1+\pi_5, 1+\pi_5^2, 1+\pi_5^3, 1+\pi_5^4, 1+\pi_5^5)\\
    &\equiv
    (1+2\pi_5+4\pi_5^2, 1+4\pi_5^2+\pi_5^3+\pi_5^4, 1+3\pi_5^3+3\pi_5^4+\pi_5^5, 1+\pi_5^4+3\pi_5^5, 1+2\pi_5^5),
\end{align*}
modulo $K_5^{\times 5}$, which implies 
\begin{equation*}
    \sigma(\beta, \gamma, \delta, \epsilon, \eta)
    \equiv (\beta^2\gamma^3\delta^4\epsilon\eta, \gamma^4\delta\eta, \delta^3\epsilon^3\eta,\epsilon\eta^3 , \eta^2 )
    \pmod{K_5^{\times 5}}.
\end{equation*}
For a prime $\fp$ not above 5, any generator $\pi_\fp$ of $\mathfrak{p}$ is not divided by $\pi_5$ so we can write 
$$\pi_\fp\equiv \zeta_4^i\beta^b\gamma^c\delta^d\epsilon^e\eta^f \pmod{\pi_5^6}.
$$
A (multiplicative) 
$\bF_5$-vector space  $\lara{\beta,\gamma,\delta,\epsilon,\eta}$ is decomposed by eigenvectors $\lcrc{\epsilon\eta^2, \gamma\delta\epsilon, \eta, \beta \gamma \epsilon, \delta \epsilon^4, \eta^3}$ of $\sigma$ such that 
\begin{equation*}
    \sigma(\epsilon\eta^2, \gamma\delta\epsilon, \eta,   \beta\gamma\epsilon, \delta\epsilon^4\eta^3) \equiv 
    (\epsilon\eta^2, (\gamma\delta\epsilon)^4 , \eta^2, (\beta\gamma\epsilon)^2, (\delta\epsilon^4\eta^3)^3) \pmod{K_5^{\times 5}}.
\end{equation*}
\indent 
(1) 
Since $\sigma(n) = n$ for all $n \in \bZ$, the class of $n$ in $K_5^\times/K_5^{\times 5}$ is a power of $\epsilon\eta^2$, which is the unique eigenvector with eigenvalue $+1$.
Note that
$$\epsilon\eta^2(1+\pi_5^6)^2(1+\pi_5^7)\equiv 1+\pi_5^4+2\pi_5^5+2\pi_5^6+\pi_5^7 \equiv 21 \pmod{\pi_5^8},
\quad \textrm{and} \quad
\zeta_4 \equiv 7 \pmod{\pi_5^8}.
$$ 
So for $i=0,1,2,3$,
\begin{equation*}
    \begin{array}{ccc}
    \zeta_4^i\epsilon\eta^2(1+\pi_5^6)^2(1+\pi_5^7) &\equiv & 21, 22, 3, 4 \pmod{25} \\
    \zeta_4^i\epsilon^2\eta^4(1+\pi_5^6)^4(1+\pi_5^7)^2 &\equiv & 16, 12, 9, 13 \pmod{25} \\
    \zeta_4^i\epsilon^3\eta(1+\pi_5^6)^6(1+\pi_5^7)^3 &\equiv & 11, 2, 14, 23 \pmod{25} \\
    \zeta_4^i\epsilon^4\eta^3(1+\pi_5^6)^8(1+\pi_5^7)^4 &\equiv & 6, 17, 19, 8 \pmod{25} \\
    \zeta_4^i &\equiv & 1, 7, 24, 18 \pmod{25}
    \end{array}
\end{equation*}
where
$(1+\pi_5^6)^2(1+\pi_5^7)$ is a $5^{\textrm{th}}$-power in $K_5^\times$.

(2)
Since $p\equiv 4 \pmod{5}$, $p$ splits into two primes. %For a prime $w$ above $p$, a uniformizer $\pi_w$ can be chosen in $\bQ(\sqrt{5})$ because a prime ideal in $\bQ(\sqrt{5})$ above $p$ is inert in $K$. 
For a generator $\pi_w$, $\sigma\pi_w\neq \pi_w$ but $\sigma^2\pi_w=\pi_w$. 
Hence the image of $\pi_w$ in $K_5^\times/K_5^{\times 5}$ is a product of a nontrivial power of the eigenvector $\gamma\delta\epsilon$ with eigenvalue $-1$ and a power of the eigenvector $\epsilon\eta^2$ with eigenvalue $+1$, say 
$$\pi_w=(\gamma\delta\epsilon)^c(\epsilon\eta^2)^e \pmod{K_5^{\times 5}}.$$
Also, $\pi_w\cdot \sigma\pi_w\equiv (\epsilon\eta^2)^{2e} \pmod{K_5^{\times 5}}$ and 
$\pi_w\cdot \sigma\pi_w\equiv p \pmod{K_5^{\times 5}}$ imply that the exponent $e$ is $0, 1, 2, 3, 4$ when $p \equiv 24, 9, 19, 4, 14 \pmod{25}$ respectively.
We also have 
\begin{align*}
    -u_F & \equiv 2+4\pi_5^2+4\pi_5^3+\pi_5^5  \equiv \zeta_4(1+2\pi_5^2+2\pi_5^3+3\pi_5^4+4\pi_5^5) \pmod{\pi_5^6}\\
                 & \equiv \zeta_4\gamma^2\delta^2\epsilon^2 \pmod{\pi_5^6}.
\end{align*}
Since $u_F$ is a fundamental unit of $\bQ(\sqrt{5})$, we note that another choice of a generator of the form $a' + b'\sqrt{5}$ for $a', b' \in \frac{1}{2}\bZ$ should be a product of power of $-1, u_F$, and $a + b\sqrt{5}$.
Let $\pi_w=a+b\sqrt{5}$ be a generator for $w$  with $a,b\in \frac{1}{2}\bZ$ and let $a\equiv 2^k \pmod{5}$ with $1\leq k\leq 4$. Since $$\frac{-1-\sqrt{5}}{2}(a+b\sqrt{5})=-\frac{a+5b}{2}-\left(\frac{a+b}{2}\right)\sqrt{5}$$ 
and $(-a-5b)/2\equiv 2a\pmod{5}$, we can find another generator 
$$\pi_w'=a'+b'\sqrt{5}=\left(-\frac{1+\sqrt{5}}{2}\right)^{5-k}\pi_w$$ 
of $w$, where $a'\equiv 2\pmod{5}$. 
We also note that every generator of $w$ is equivalent to one of $\pi_w'$ up to $K^{\times 5}$.

Now assume $a\equiv 2\pmod{5}$.
Then
\begin{align*}
    \zeta_4^{3} \cdot (a + b \sqrt{5}) &= (3 + 2\pi_5^4 + 4\pi_5^5 + O(\pi_5^6))(a + b(2 \pi_5^2 + 2 \pi_5^3 +  \pi_5^4 + O(\pi_5^6))) \\
    &= 1+ b\pi_5^2+O(\pi_5^3)
\end{align*}
implies that $\pi_w=(\gamma\delta\epsilon)^b(\epsilon\eta^2)^e$ in $K^{\times}/K^{\times 5}$. 
% $b$ determines the exponent of $(\gamma \delta \epsilon)$, and $e$ is already determined by $p$ mod 25.
This induces the first row of the table. 
The other rows are determined by the relation between $\pi_w'$ and $\pi_w$ and the value of $-(1+\sqrt{5})/2$ in $K_5^\times/K_5^{\times 5}$.
% If $a\equiv 2^k\pmod{5}$ with $k\geq 2$, then an induction on $k$ by
% $$a+b\sqrt{5}=\frac{-1-\sqrt{5}}{2}\times \left(\frac{a-5b}{2}+\frac{b-a}{2}\sqrt{5}\right)$$
% implies that $\pi_w$ is equivalent to $\displaystyle{(\gamma\delta\epsilon)^2(2^{k-1}+(3b-2^{k-1})\sqrt{5})=(\gamma\delta\epsilon)^{3^{k-1}b}(\epsilon\eta^2)^e}$ in $K^{\times}/K^{\times 5}$.
\end{proof}

In the next section, we will need the images of $\lcrc{\zeta_5, 1 \pm \zeta_5, 2}$ in $K_\fp^\times/K_\fp^{\times 5}$ also.
We begin with $\fp = 2$.
Recall that $K_2^\times / K_2^{\times 5}\cong \lara{2, \zeta_5} = \lara{\alpha, \beta}$ in (\ref{eqn: local field not 5}).
\begin{lemma} \label{lem:localover2}
(1) The image of $(\zeta_5, 1 + \zeta_5, 1-\zeta_5, 2)$ in $K_2^\times/K_2^{\times 5}$ is $(\beta, \beta^3, \beta^3, \alpha)$. \\
(2) The images of odd integers and prime elements $\pi_w = a+ bu_F$ for $a, b \in \bZ$ in $K_2^\times/K_2^{\times 5}$ are trivial.
\end{lemma}
\begin{proof}
(1)
To describe $2$-expansions of elements of $K_2$, we fix an isomorphism
\begin{equation*}
    \bF_{16} \cong \bF_2[t]/(t^4 + t + 1).
\end{equation*}

We choose an embedding of $K$ in $K_2$ which sends $\zeta_5 \in K$ to $t^3 \in \bF_{16}$.
Since
\begin{equation*}
    (t^3+1)(t^2+t+1) = t^3 + t, \qquad
    (t^2+t+1)^3 = 1, \qquad t^9 = t^3 + t,
\end{equation*}
we know that $(1 + \zeta_5) \zeta_3 = \zeta_5^3$ in $K_2$. Since $\zeta_3$ is trivial in $K_2^{\times}/K_2^{\times 5}$, the image of $(1 + \zeta_5)$ in $K_2^{\times}/K_2^{\times 5}$ is $\beta^3$.
Also, the 2-expansion of the image of $(1-\zeta_5)$ in $K_2$ is
\begin{equation*}
    1 - \zeta_5 = 1 + t^3(1 + 2 + O(2^2)) = (1+t^3)(1 + (1+t^3)^{-1}t^3 2 + O(2^2)).
\end{equation*}
Hence the image of $(1 - \zeta_5)$ in $K_2^\times/K_2^{\times 5}$ is $\beta^3$ also.

(2) Since $U^{(1)}$ vanishes in $K_2^\times/K_2^{\times 5}$, every odd integer maps to the trivial element in $K_2^\times/K_2^{\times 5}$. In $K_2$, one has
\begin{equation*}
    \sqrt{5} = 1 + (t^2 + t)2 + O(2^2) \quad \textrm{and} \quad
    u_F = (t^2 + t + 1) + O(2).
\end{equation*}
Therefore, the image of $a + bu_F$ in $\bF_{16}^{\times}$ is contained in $\lcrc{t^2+t+1, t^2+t, 1}$ which is the group generated by $\zeta_3$.
\end{proof}

\begin{lemma} \label{lem: local over not5}
Let $p \neq 2$ be a rational prime inert in $K/\bQ$ and let $\pi_w$ be a prime element defined by $a + b \sqrt{5}$ for $a, b \in \frac{1}{2}\bZ$. \\
(1) For $\fp = (p)$ or $(\pi_w)$, the image of $\lcrc{\zeta_5, 1+ \zeta_5, 1-\zeta_5}$ in $K_\fp^\times/K_\fp^{\times 5}$ is in $\lara{\beta_\fp}$. \\
(2) For $\fp = (p)$, the images of rational primes relatively prime to $\fp$ and prime elements $\pi_{w'} = a'+ b'\sqrt{5}$ for $a', b' \in \frac{1}{2} \bZ$  are trivial in $K_\fp^\times/K_\fp^{\times 5}$. \\
(3) For $\fp = (\pi_w)$, the images of rational primes relatively prime to $\fp$ and a prime element $\pi_{\overline{w}} := a - b\sqrt{5}$ are trivial in $K_\fp^\times/K_\fp^{\times 5}$.
\end{lemma}
\begin{proof}
(1) We recall that $K_p^\times \cong p^{\bZ} \times \mu_{p^4-1} \times U^{(1)}$ and
$K_w^\times \cong \pi_w^{\bZ} \times \mu_{p^2 - 1} \times U^{(1)}$, i.e. $K_\fp^\times/K_\fp^{\times 5} = \lara{\alpha_\fp, \beta_\fp}$ for $\fp = (p)$ or $(w)$ in (\ref{eqn: local field not 5}).
Especially, $U^{(1)}$-part vanishes in $K_\fp^\times/K_\fp^{\times 5}$. 
Since $\zeta_5, 1 \pm  \zeta_5$ are not divided by $\fp$, their images are in $\lara{\beta_\fp}$.

(2) Every rational integer relatively prime to $p$ and $\pi_{w'}$ maps to $\bF_{p^2}^\times$ modulo $p$. Since the fifth-power map on $\bF_{p^2}^\times$ is bijective, every element maps to $\bF_{p^2}^\times$ vanish in $K_p^\times/K_p^{\times 5}$.

(3) Similarly, every integer and $\pi_{\overline{w}}$ maps to $\bF_{p_w}^\times$ where $p_w$ is the rational prime divided by $\pi_w$.
\end{proof}

% \iffalse
% \begin{lemma}
% Let $v$ be a prime above a rational prime $p$ equivalent to 1 modulo 5.
% In $K_v^\times/K_v^{\times 5} \cong \lara{\pi_v, \zeta_{5^n}} =: \lara{\alpha, \beta}$, \\
% \textcolor{blue}{
% (1) The image of $\lcrc{1 + \zeta_5, 1-\zeta_5}$, that of 2, ... ... \\
% (2) The image of $\pi_v'$ which is lying over same prime with $\pi_v$ is in $\lara{\beta}$. \\
% (3) The image of odd integers and generators $\pi_w = a+b\sqrt{5}$ of primes above a rational prime equivalent to 4 modulo 5 is  \\
% (4) The image of other split prime is }
% \end{lemma}
% \begin{proof}
% We note that $n$ is the integer satisfying $5^n \| (p-1)$.

% (1) 
% An image of 2 is in $\lara{\beta}$ and is trivial if and only if $2^{(p-1)/5} \equiv 1$. (For example, $p =151$ satisfies this condition.)
% Similarly, an image of $\zeta_5$ in $K_v^\times/K_v^{\times 5}$ is in $\lara{\beta}$ and trivial if and only if $5^2 \mid (p-1)$.

% (2) 
% \end{proof}
% \fi

\subsection{The root numbers}

We recall the result of  \cite{Sto2} on the root numbers of $y^2 = x^l + A$.
\begin{theorem}[{\cite[Theorem 3.2]{Sto2}}]
The root number $w(A)$ of the curve $y^2 = x^l + A$ over $\bQ$ where $A$ is a $2l$-th power free integer not divisible by $l$, is given by
\begin{equation*}
w(A) = \left\{ \begin{array}{cc}
\lbrb{\frac{2A v_A}{l}} & \textrm{if $l \mid q_l(A)$}, \\
-\lbrb{ \frac{2 q_l(A) v_A}{l} } & \textrm{if $l \nmid q_l(A)$},
\end{array} \right.
\end{equation*}
where $q_l(A) = (A^{l-1} -1 )/ l$ and $v_A = 2^{f_2(A)} \prod_{p \mid A, p \neq 2} p$ where $f_2$ is given by
\begin{equation*}
f_2(A) = \left\{ 
\begin{array}{lll}
0 & \textrm{if $e = 2l-2$ and $B \equiv 1 \pmod{4}$},\\
1 & \textrm{if $e < 2l-2$ and is even and $B \equiv 1 \pmod{4}$}, \\
2 & \textrm{if $e$ is even and $B \equiv -1 \pmod{4}$}, \\
3  & \textrm{if $e$ is odd.}
\end{array}
\right.
\end{equation*}
for $A = 2^eB$ with $B$ odd.
\end{theorem}

In this paper, we only need the following special case.

\begin{corollary} \label{cor:rootnumber}
For an odd square-free integer $m$, the root number $w(m^2)$ of the hyperelliptic curve $y^2 = x^5 + m^2$ over $\bQ$ is given by
\begin{equation*}
    w(m^2) = \left\{ \begin{array}{cc}
        +1 & \textrm{if } m \equiv 
        1, 2, 4, 6, 12, 13, 19, 21, 23, 24 \pmod{25}, \\
        -1 & \textrm{if } m \equiv 
        3, 7, 8, 9, 11, 14, 16, 17, 18, 22 \pmod{25}.
        \end{array} \right.
\end{equation*}
\end{corollary}

\section{Descent for Jacobian of hyperelliptic curves} \label{sec:descent}

We recall the general facts on the descent for Jacobian of hyperelliptic curves of odd prime degree. The main reference is \cite{Sch}.

Let $p$ be an odd prime, let $K$ be a number field containing $\zeta_p$, and let $C$ be a curve defined by an equation $y^p = f(x)$.
Let $J$ be the Jacobian of $C$ and consider an endomorphism $\phi$ of $J$. 
The $\phi$-Selmer group of $J/K$ is defined by
\begin{equation*}
    \Sel_{\phi}(J/K) := 
    \ker \lbrb{ H^1(K, J[\phi]) \to \prod_\fp H^1(K_\fp , J) }.
\end{equation*}
where $\fp$ is taken over all primes of $K$.
Following the Schaefer's idea, instead of using the first cohomology group we will use more concrete object which we will describe as follows.
Assume that $J[\phi]$ has a prime power exponent $q$.
We define
\begin{equation*}
    L := K[T]/(f(T)), \qquad
    H := \ker \lbrb{ \Norm : L^\times/L^{\times q} \to K^\times/K^{\times q} }.
\end{equation*}
Let $\lambda : J \to \widehat{J}$ be the canonical polarization of $J$ and let $\widehat{\phi}$ be the dual isogeny of $\phi$.
Let $\Psi := \lambda^{-1}(\widehat{J}[\widehat{\phi}]) \subset J[q]$ and choose a $G_K$-invariant set of divisor classes that generate $\Psi$.
We also define $\Div^0_{\perp}(C)$ as a set of degree zero divisors of $C$ with support not intersecting with the generating set of $\Psi$.
For each element of $J(K)$, we may choose its representative in $\Div^0_{\perp}(C)$.
There is a map 
\begin{equation*}
    F : \Div^0_{\perp}(C) \to L^\times
\end{equation*}
which induces $F : J(K)/\phi J(K) \to L^\times/L^{\times q}$ by \cite[Lemma 2.1, Theorem 2.3]{Sch}.

Now we consider our cases $p = 5$, $K = \bQ(\zeta_5)$, $C_{m^2} : y^2 = x^5 + m^2$ and $\phi = (1 - \zeta_5)$ where $\zeta_5(x_0, y_0) := (\zeta_5 x_0, y_0)$.
We note that the class number of $K$ is one and there is a fundamental unit $(1 + \zeta_5)$.
Let $J_{m^2}$ be the Jacobian of $C_{m^2}$.
The polynomial $f(T) = T^2 - m^2$ is reducible so we have $L \cong K \oplus K$, and the norm map is given by $(k_1, k_2) \to k_1k_2$.
After identifying $H$ with $K^\times$, we have
\begin{equation*}
    H^1(K, J_{m^2}[\phi] ; S) \cong K(S, 5)
\end{equation*}
where $K(S, 5)$ is a subset of $K^\times/K^{\times 5}$ consisting of elements trivial outside $S$.
Since the set of bad primes $S$ consists of the primes above $10m$, we note that $K(S, 5)$ is generated by
\begin{equation*}
    \zeta_5, 1 + \zeta_5, 2, 1-\zeta_5
\end{equation*}
and prime elements dividing $m$.
We also have $\lambda^{-1}(\widehat{J_{m^2}}[\widehat{\phi}]) = J_{m^2}[\phi]$
and $\lcrc{(0,\pm m) - \infty }$ forms a basis of $J_{m^2}[\phi]$ by \cite[Proposition 3.1, 3.2]{Sch}.
Furthermore, we have %$H \cong H^1(K, J_{m^2}[\phi])$ and
\begin{equation} \label{eqn:Sel Sch}
    \Sel_\phi(J/K) \cong \bigcap_{\fp \in S} (i_\fp^{-1}\circ F_\fp) \lbrb{J_{m^2}(K_\fp)/\phi J_{m^2}(K_\fp) },
\end{equation}
where $i_\fp$ is a natural map $L^\times \to L_\fp^\times$, by \cite[Proposition 3.4]{Sch}.
For the concrete computation, we remind that 
\begin{equation} \label{eqn:local dimension}
    \dim_{\bF_p} (J_{m^2}(K_\fp)/\phi J_{m^2}(K_\fp)) = \left\{
    \begin{array}{cc}
        3 & \textrm{if } \fp \mid 5,  \\
        1 & \textrm{otherwise},
    \end{array}
    \right.
\end{equation}
by \cite[Corollary 3.6]{Sch}.
This result guides us when we stop finding the independent points of $J_{m^2}(K_\fp)/\phi J_{m^2}(K_\fp)$.
Also, for $D = Q_1 + \cdots + Q_r - r\infty$ where $Q_i$ are $K$-conjugates with $x(Q_i) \neq 0$,
\begin{equation*}
    F_\fp([D]) = \prod_{i = 1}^r (y(Q_i) - m)
\end{equation*}
and for $D = (0, \pm m) - \infty$, 
\begin{equation*}
    F_\fp([D]) = \pm 2m
\end{equation*}
by \cite[Proposition 3.3]{Sch}.
We remark that 
\begin{equation*}
    \rk(J_{m^2}(\bQ)) = \dim_{\bF_5}(J_{m^2}(K)/\phi J_{m^2}(K)) - \dim_{\bF_5} J_{m^2}(K)[\phi],
\end{equation*}
by \cite[Corollary 3.7, Proposition 3.8]{Sch}.

One of the main goals of the paper is computing the Selmer group of Jacobian of $C_{m^2}$.
\begin{proposition} \label{prop:map F on 5}
Let $m$ be an odd integer and let $J_{m^2}$ be a Jacobian of $C_{m^2}$.
Under the identifications of $K_\fp^\times/K_\fp^{\times 5}$ as in (\ref{eqn:local K5}) and (\ref{eqn: local field not 5}), we have
\begin{align*}
    F_5(J_{m^2}(K_5)/\phi J_{m^2}(K_5)) &=
        \lara{\delta, \epsilon, \eta} \qquad \textrm{if } m \equiv \pm 1, \pm 7 \pmod{25}
\end{align*}
and $F_\fp(J_{m^2}(K_\fp)/ \phi J_{m^2}(K_\fp)) = \lara{\alpha_\fp}$ for all other $\fp$.% We note that $\alpha$ depends on the prime $\fp$.
\end{proposition}
\begin{proof}
In the proof, we denote $J$ by $J_{m^2}$.
The $F_5$-case is a generalization of \cite[Proposition 3.12]{Sch}.
We recall that
\begin{equation*}
K_5^\times/K_5^{\times 5} \cong \lara{\pi_5, 1+\pi_5, 1+ \pi_5^2, 1+ \pi_5^3, 1+ \pi_5^4, 1+ \pi_5^5} :=
\lara{\alpha, \beta, \gamma, \delta, \epsilon, \eta}
\end{equation*}
and every element of $K_5^\times$ which is one modulo $\pi_5^6$ is fifth power.
When $m^2 \pm 1 \equiv 0 \pmod{25}$, either $y^2 - m^2 \equiv 1 \pmod{\pi_5^6}$ or $m^2 - y^2 \equiv 1 \pmod{\pi_5^6}$ has solutions $\pi_5^i$ for $i = 3, 4, 5$.
Hence, in each case, there is an $x_i$ such that $[(x_i, \pi_5^i) - \infty ]$ for $i = 3, 4, 5$ is the point of $J(K_5)/\phi J(K_5)$.
The value of $F_5((x_i, \pi_5^i) - \infty)$ is determined by the image of $\pi_5^i + m$ in $K_5^\times/K_5^{\times 5}$.
For $m \equiv \pm 1, \pm 7 \pmod{25}$, the images of $\pi_5^i+m$ in $U^{(2)}$ are
\begin{equation*}
    (1 + \pi_5^i), \qquad (1 - \pi_5^i), \qquad 
    \zeta_4^3(7 + \pi_5^i), \qquad \zeta_4^3(7 - \pi_5^i)
\end{equation*}
respectively.
Computing the $\pi_5$-expansion, we get
\begin{equation*}
    \begin{array}{cccccc}
         & y + 1 & y - 1 & y + 7 & y - 7  \\
        {[}(x_3, \pi_5^3) - \infty{]} & \delta  & \delta^{-1} & \delta^3 & \delta^2 \\
        {[}(x_4, \pi_5^4) - \infty{]} & \epsilon & \epsilon^{-1} & \epsilon^3 & \epsilon^2 \\
        {[}(x_5, \pi_5^5) - \infty{]} & \eta & \eta^{-1} & \eta^3 & \eta^2
    \end{array}
\end{equation*}
Together with (\ref{eqn:local dimension}) we have
\begin{equation*}
    F_5(J(K_5)/\phi J(K_5)) = \lara{\delta, \epsilon, \eta}.
\end{equation*}

Again by (\ref{eqn:local dimension}) for $\fp \nmid 5$, we have $\dim_{\bF_5}(J(K_\fp)/\phi J(K_\fp)) = 1$.
By Lemma \ref{lem:localover2}, arbitrary odd integer $m$ maps to 1 in $K_2^\times/K_2^{\times 5} \cong \lara{2, \zeta_{5}} =  \lara{\alpha_2, \beta_2}$. Hence, 
\begin{equation*}
    \begin{array}{cccccc}
         & y + m & y - m  \\
        {[}(0, m) - \infty{]} & 2  & 2^{-1}
    \end{array}
\end{equation*}
and $F_2(J(K_2)/\phi J(K_2))$ is $\lara{\alpha_2}$.
Similarly for $\fp \nmid 10$, the image of 2 in $K_\fp^\times/K_\fp^{\times 5}$ is trivial by Lemma \ref{lem: local over not5}.
So
\begin{equation*}
    \begin{array}{cccccc}
         & y + m & y - m  \\
        {[}(0, m) - \infty{]} & m  & m^{-1}
    \end{array}
\end{equation*}
shows that $F_\fp(J(K_\fp)/\phi J(K_\fp)) = \lara{\alpha_\fp}$.
\end{proof}

\begin{remark} \label{rem:}
We note that Proposition \ref{prop:map F on 5} is enough to prove the main theorem, but the same strategy gives $F_5(J_{m^2}(K_5)/\phi J_{m^2}(K_5))$ when one knows the generators of $J_{m^2}(K_5)/\phi J_{m^2}(K_5)$.
For example, 
\begin{equation*}
    (-\pi_5, 2 + 3\pi_5^4+2\pi_5^5), \qquad
    (1 , \pi_5^2 + \pi_5^3 + 3\pi_5^4), \qquad
    (2,1)
\end{equation*}
are solutions of $y^2 \equiv x^5 + m^2 \pmod{\pi_5^6}$ when $m \equiv \pm 12 \pmod{25}$.
Therefore,
\begin{align*}
    (\zeta_4^2(2 + 3\pi_5^4+2\pi_5^5 & + 12), 
    \zeta_4^3(\pi_5^2 + \pi_5^3 + 3\pi_5^4 + 12),
    \zeta_4(1 + 12) )  \\
    & \equiv
    (1 + 4\pi_5^5, 1 + 3\pi_5^2 + 3\pi_5^3 + \pi_5^4 +4 \pi_5^5 , 1+ 2\pi_5^4 + 4\pi_5^5 ) \pmod{\pi_5^6} \\
    & \equiv (\eta^4, \gamma^3 \delta^3 \epsilon , \epsilon^2\eta^4) \quad \textrm{ in } K_5^\times/K_5^{\times 5}.
\end{align*}
Hence,
\begin{equation*}
    F_5(J_{m^2}(K_5)/\phi J_{m^2}(K_5)) = \lara{\gamma \delta, \epsilon, \eta}
\end{equation*}
when $m \equiv \pm 12 \pmod{25}$.
Similarly we can compute $F_5(J_{m^2}(K_5)/\phi J_{m^2}(K_5))$ for other cases.
Also, Lemma \ref{lem:localover2} and \ref{lem: local over not5} describe an image of prime element not lying above $p \equiv 1 \pmod{5}$.
Therefore, we can calculate the Selmer group of $J_{m^2}$ when $m$ is square-free and
\begin{itemize}
    \item[(1)] if $p$ divides $m$ then $p \not\equiv 1 \pmod{5}$,
    \item[(2)] there is at most one prime divisors $p$ of $m$ such that $p \equiv 4 \pmod{5}$,
\end{itemize}
even though we do not fully describe the result.
We will give an example in the end of this section.
\end{remark}

\begin{proposition} \label{prop:descent not over 5}
    Let $m$ be a square-free integer satisfying the above two conditions (1), (2) and let $\fp \nmid 5$ be a prime of $K$ dividing $m$. 
    %We choose a generator of each prime dividing $m$ as in Lemma \ref{lem: local over not5}.
    % \begin{align*}
    %     \pi = \left\{ \begin{array}{cc}
    %         p & \textrm{if $(\pi) = (p)$ is inert in $K/\bQ$ } \\
    %         a + b\sqrt{5} & \textrm{if $(\Norm^K_{\bQ}(\pi))$ is split into two primes}
    %     \end{array} \right.
    % \end{align*}
    % $p$ when the prime is inert in $K/\bQ$ and $a + b\sqrt{5}$ for $a, b \in \frac{1}{2}\bZ$ when the prime splits by two primes.
    Then,
    $(i_\fp^{-1} \circ F_\fp)(J_{m^2}(K_\fp)/\phi J_{m^2}(K_\fp))$ contains $2$ and prime generators dividing $m$ chosen as in Lemma \ref{lem: local over not5}.
\end{proposition}
\begin{proof}
This is a direct consequence of Lemma \ref{lem: local over not5} and Proposition \ref{prop:map F on 5}.
\end{proof}

\begin{corollary} \label{cor: phiSelmer}
For a rational prime $p$ and the Jacobians $J_{p^2}$, we have
\begin{equation*}
    \dim_{\bF_5} \Sel_{\phi} (J_{p^2}/\bQ) = 2, \qquad \textrm{if } p \equiv 7, 8 \pmod{25}.
\end{equation*}
When $p \equiv 24 \pmod{25}$, there is a generator $\pi_w$ of $w$ above $p$ satisfies $\pi_w = a + b\sqrt{5}$ for $a, b \in \frac{1}{2}\bZ$. Then,
\begin{equation*}
    \dim_{\bF_5}\Sel_{\phi}(J_{p^2}/\bQ) = \left\{
    \begin{array}{cc}
        1 & b \not\equiv 0 \pmod{5}, \\
        3 & b \equiv 0 \pmod{5}.
    \end{array}
    \right.
\end{equation*}
\end{corollary}
\begin{proof}
In the proof, we denote $J$ by $J_{p^2}$.
We first consider the case of $p \equiv 7, 8 \pmod{25}$.
We recall that $i_5 : K(S, 5) \to K^\times/K^{\times 5}$, and $K(S, 5)$ is generated by $\zeta_5, 1+\zeta_5, 2, 1-\zeta_5$ and a prime $p$ which is inert in $K/\bQ$.
Since
\begin{equation*}
    i_5(\zeta_5, 1+\zeta_5, 2, 1-\zeta_5, 7,8) = 
    (\beta \gamma \epsilon, \beta^2\gamma^4 \delta^2 \epsilon^4, \epsilon^3 \eta, \alpha, 1, \epsilon^4\eta^3)
\end{equation*}
by Lemma \ref{lem:localover5}, we have
\begin{equation*}
    F_5(J(K_5)/ \phi J(K_5)) = \lara{\delta, \epsilon, \eta}, \qquad
    \im i_5 = \lara{\beta\gamma\epsilon, \beta^2\gamma^4\delta^2\epsilon^4, \epsilon^3\eta, \alpha},
\end{equation*}
together with Proposition \ref{prop:map F on 5}.
A sort of linear algebra 
\iffalse
(More concretely, since 
\begin{equation*}
    \lara{\beta \gamma \epsilon, \beta^2 \gamma^4 \delta^2 \epsilon^4, \epsilon^3\eta} \cap \lara{\delta, \epsilon, \eta} \supset \lara{\epsilon^3 \eta},
    \qquad
    \lara{\beta \gamma \epsilon, \beta^2 \gamma^4 \delta^2 \epsilon^4, \epsilon^3\eta} + \lara{\delta, \epsilon, \eta} = \lara{\beta, \gamma, \delta, \epsilon, \eta}
\end{equation*}
we have $\dim \lara{\beta \gamma \epsilon, \beta^2 \gamma^4 \delta^2 \epsilon^4, \epsilon^3\eta} = 3$ and the dimension of the intersection is one) \fi
shows that
\begin{equation*}
    \im i_5 \cap F_5(J(K_5)/\phi J(K_5))
    = \lara{\epsilon^3\eta},
\end{equation*}
and
\begin{equation*}
    (i_5^{-1}\circ F_5)(J(K_5)/\phi J(K_5)) = \lara{2, p}.
\end{equation*}
By Proposition \ref{prop:map F on 5}, $F_\fp(J(K_\fp)/\phi J(K_\fp)) = \lara{\alpha_{\fp}}$ for a prime $\fp$ not above 5.
Now, Proposition \ref{prop:descent not over 5} gives 
\begin{equation*}
    (i_2^{-1} \circ F_2)(J(K_2)/\phi J(K_2)) \supset \lara{2, p}, \qquad
    (i_p^{-1} \circ F_p)(J(K_p)/\phi J(K_p)) \supset \lara{2, p},
\end{equation*}
which shows that $\dim_{\bF_5}\Sel_\phi(J/\bQ) = 2.$

When $p \equiv 24 \pmod{25}$, we choose the generators $\pi_w, \pi_{\overline{w}}$ above $p$ by $a \pm b\sqrt{5}$ for $a, b \in \frac{1}{2}\bZ$.
We still have $F_5(J(K_5)/ \phi J(K_5)) \cong \lara{\delta, \epsilon, \eta}$. 
By Lemma \ref{lem:localover5}, the images under $i_5$ of the generators above $p \equiv 24$ are in  $\lara{\gamma \delta \epsilon}$ and trivial when $b \equiv 0 \pmod{5}$.
Hence,
\begin{equation*}
    \im i_5 \subset \lara{\beta \gamma \epsilon, \beta^2 \gamma^4 \delta^2 \epsilon^4, \epsilon^3 \eta, \alpha, \gamma \delta \epsilon}.
\end{equation*}
Since $(\beta\gamma\epsilon)^3(\beta^2\gamma^4\delta^2\epsilon^4)(\gamma \delta \epsilon)^2$ is trivial, the dimension of the space in the right hand side is 4.
Hence, the similar argument gives
\begin{equation*}
    \im i_5 \cap F_5(J(K_5)/\phi J(K_5))
    = \lara{\epsilon^3\eta},
\end{equation*}
and
\begin{equation*}
    (i_5^{-1}\circ F_5)(J(K_5)/\phi J(K_5)) = 
    \left\{
    \begin{array}{ll}
    \lara{2}     &  \textrm{if } b\not\equiv 0  \pmod{5}, \\
    \lara{2, \pi_w, \pi_{\overline{w}}} &  \textrm{if } b \equiv 0 \pmod{5}.
    \end{array}
    \right.
\end{equation*}
Together with Proposotion \ref{prop:descent not over 5}, we know that the dimension of the Selmer group $\Sel_\phi(J_{p^2}/\bQ)$ is 1 or 3, and dimension 3 if and only if $b\equiv 0 \pmod{5}$.
\end{proof}

\begin{proof}[Proof of Theorem \ref{thm:main MW group}]
By the Dirichlet theorem on arithmetic progressions for number fields, there are infinitely many primes in a ray class modulo an ideal.
Let us denote two real embeddings by $\sigma_1, \sigma_2$.
For a modulus $(50) \cdot \sigma_1\sigma_2$ and a ray class $(2+ \sqrt{5})$, there are infinitely many prime elements $\pi$ which are congruent modulo $(50)\cdot \sigma_1\sigma_2$ to one of $u_F^{2n}(2 + \sqrt{5})$ where $u_F = (1+\sqrt{5})/2$.

Using an integral basis $\lcrc{1, u_F}$ of $\cO_{F}$, we may write
\begin{align*}
    \pi =  u_F^{2n}(2 + \sqrt{5}) + 50z_1 + 50z_2 u_F
\end{align*}
for some $z_1, z_2 \in \bZ$. 
Then, the norm of $\pi$ is $-1 \pmod{25}$.
Let $a_{n}$ and $b_n$ be integers satisfying
\begin{align*}
    u_F^n = a_n + b_nu_F.
\end{align*}
Then, 
\begin{align*}
    \pi &= u_F^{2n}\lbrb{2 + \sqrt{5} \pm 50z_1\lbrb{a_{-2n}+b_{-2n}u_F} 
    \pm 50z_2\lbrb{ a_{-2n+1} +b_{-2n+1}u_F }} \\
    &= u_F^{2n} (
    2 +\sqrt{5} \pm 25 (z_1 (2a_{-2n} + b_{-2n}) + z_2 (2a_{-2n+1} +b_{-2n+1})
     + \sqrt{5}(z_1 b_{-2n} + z_2b_{-2n+1} ))).
\end{align*}
For a rational prime $p \equiv 24 \pmod{25}$ divided by $\pi$, there is a generator of $(\pi)$ satisfying the condition of Corollary \ref{cor: phiSelmer} with $b \not \equiv 0 \pmod{5}$.
From the exact sequence
\begin{align*}
    \xymatrix{0 \ar[r] & \displaystyle{ \frac{J_{p^2}(\bQ)}{\phi J_{p^2}(\bQ)} } \ar[r] & \Sel_{\phi}(J_{p^2}/\bQ) \ar[r] & \Sha(J_{p^2}/\bQ)[\phi] \ar[r] & 0 }
\end{align*}
and Corollary \ref{cor: phiSelmer}, one can deduce that $J_{p^2}(\bQ) \cong \bZ/5\bZ$.

Also, for a prime $p \equiv 7, 8 \pmod{25}$ we have
\begin{equation*}
    \bZ/5\bZ \leq J_{p^2}(\bQ) \leq \bZ/5\bZ \times \bZ, \qquad w(p^2) = -1
\end{equation*}
by Corollary \ref{cor: phiSelmer} and Corollary \ref{cor:rootnumber}.
This proves the second part of the theorem.
\end{proof}

We note that the machinery also works for the totally split primes, even though one need to compute everything directly.
\begin{proposition} \label{prop:101}
The Mordell--Weil rank of $J_{101^{2}}/\bQ$ is zero.
\end{proposition}
\begin{proof}
We will show that $\dim_{\bF_5}\Sel_{\phi}(J_{101^2}/\bQ) = 1$.
We note that Sagemath \cite{Sag} runs most of computation in the proof.
Let $\fp_j$ for $j = 1, 2, 3, 4$ be a prime ideal of $K$ above $p = 101$, and let us choose generators $\pi_j$ by
\begin{align*}
    \zeta_5^3 + 3\zeta_5^2 - \zeta_5 + 1, \qquad
    3\zeta_5^3 + 4\zeta_5^2 + 2\zeta_5 + 2, \qquad
    -4\zeta_5^3 - 2\zeta_5^2 - \zeta_5 - 2, \qquad
    -2\zeta_5^3 - \zeta_5^2 + 2\zeta_5.
\end{align*}
We note that $\pi_1\pi_2\pi_3\pi_4 = 101$.
Also, 
\begin{align*}
    K(S, 5) = \lara{2, \zeta_5, 1 + \zeta_5, 1-\zeta_5, \pi_1, \pi_2, \pi_3, \pi_4}.
\end{align*}
Now we want to compute the image of $i_{1} := i_{\pi_1} : K(S, 5) \to K_{\fp_1}^\times/K_{\fp_1}^{\times 5}$ of the above generators.
In section \ref{sec:Preliminaries} we showed that $K_{\fp_1}^\times/K_{\fp_1}^{\times 5}$ is generated by two elements $\alpha_{\fp_1}, \beta_{\fp_1}$ which is $\pi_{\fp_1}$ and $\zeta_{25}$, respectively.
Let  $\rho_1 : \cO_{K, \fp_1} \to \cO_{K, \fp_1}/\fp_1 \cO_{K, \fp_1} \cong \bF_{101}$ be a projection map. Then,
\begin{align*}
    \rho_1( 2, \zeta_5, 1 + \zeta_5, 1-\zeta_5, \pi_2, \pi_3, \pi_4) = (2, 95, 96, 7, 92, 89, 81).
\end{align*}
We also denote $\rho_1$ as a composition of the previous map and the quotient $\bF_{101}^\times \to \bF_{101}/\bF_{101}^{\times 5}$. Then, we know that 
\begin{align*}
    \rho_1(2, \zeta_5, 1 + \zeta_5, 1-\zeta_5, \pi_2, \pi_3, \pi_4) = (\overline{2}, \overline{1}, \overline{3}, \overline{3}, \overline{8}, \overline{2}, \overline{2}).
\end{align*}
Note that $\overline{2}^3 = \overline{8}$ and $\overline{2}$ is a multiplicative inverse of $\overline{3}$.
Since the elements above are not divided by $\pi_1$, we can describe the images of elements in $K(S, 5)$ in $K_{\fp_1}^\times/K_{\fp_1}^{\times 5}$.
Now 
\begin{equation*}
    \begin{array}{cccccc}
         & y + m & y - m  \\
        {[}(0, m) - \infty{]} & 2m  & (2m)^{-1}
    \end{array}
\end{equation*}
Therefore, $F_{\fp_1}(J(K_{\fp_1}) / \phi J(K_{\fp_1}))$
is generated by the product of $\alpha_{\fp_1}$ and the image of 2. Hence, 
\begin{align*}
    (i_1^{-1}\circ F_{\fp_1} )(J(K_{\fp_1}) / \phi J(K_{\fp_1})) =
    \lara{ 2\pi_1, \zeta_5, 2(1+\zeta_5), 2(1-\zeta_5), 2^2\pi_2, 2^4\pi_3, 2^4\pi_4}.
\end{align*}
Similarly, we have
\begin{align*}
    \rho_2(2, \zeta_5, 1 + \zeta_5, 1- \zeta_5, \pi_1, \pi_3, \pi_4)
    = (\overline{2}, \overline{1}, \overline{3}, \overline{8}, \overline{2}, \overline{8}, \overline{2}),
\end{align*}
so $F_{\fp_2}(J(K_{\fp_2}) / \phi J(K_{\fp_2}))$ is generated by the product of $\alpha_{\fp_2}$ and the image of 2.
Hence,
\begin{align*}
    (i_2^{-1}\circ F_{\fp_2} )(J(K_{\fp_2}) / \phi J(K_{\fp_2})) = \lara{2\pi_2, \zeta_5, 2(1+\zeta_5), 2^2(1-\zeta_5), 2^4 \pi_1, 2^2 \pi_3, 2^4\pi_4}.
\end{align*}
Also, 
\begin{align*}
    \rho_3(2, \zeta_5, 1 + \zeta_5, 1- \zeta_5, \pi_1, \pi_2, \pi_4)
    &= (\overline{2}, \overline{1}, \overline{3}, \overline{3}, \overline{2}, \overline{2}, \overline{8}) \\
    \rho_4(2, \zeta_5, 1 + \zeta_5, 1- \zeta_5, \pi_1, \pi_2, \pi_3)
    &= (\overline{2}, \overline{1}, \overline{2}, \overline{8}, \overline{8}, \overline{2}, \overline{2})
\end{align*}
and
\begin{align*}
    (i_3^{-1}\circ F_{\fp_3} )(J(K_{\fp_3}) / \phi J(K_{\fp_3})) &=
    \lara{2 \pi_3, \zeta_5, 2(1+\zeta_5), 2(1-\zeta_5), 2^4 \pi_1, 2^4\pi_2, 2^2 \pi_4}, \\
    (i_4^{-1}\circ F_{\fp_4} )(J(K_{\fp_4}) / \phi J(K_{\fp_4})) &=
    \lara{2\pi_4, \zeta_5, 2^4(1+\zeta_5), 2^2(1-\zeta_5), 2^2 \pi_1, 2^4\pi_2, 2^4\pi_3}.
\end{align*}
We denote each vector space $(i_j^{-1} \circ F_{\fp_j})(J(K_{\fp_j}) /\phi J(K_{\fp_j}))$ over $\bF_5$ by $V_j$ for $j = 1, 2, 3, 4$.
% Since $2\pi_1 \cdot 2\pi_1^{-1} \in V_1 + V_2$, we know that $V_1 + V_2$ is the whole space $K(S, 5)$ of dimension 8. Since $\dim V_1 = \dim V_2 = 7$, we also know that $\dim (V_1 \cap V_2) = 6$.
% We can check that
% \begin{align*}
%     V_1 \cap V_2 = \lara{\zeta_5, 2(1+\zeta_5), 2\pi_4^4, 2^2 \pi_1 \pi_2 \pi_3, 2^3(1-\zeta_5)\pi_2, 2(1-\zeta_5)^2\pi_3}.
% \end{align*}
% We briefly explain how we get it. The first three are easy to find. Eliminating them, we need to find basis of 
% \begin{align*}
%     \lara{2 \pi_1, 2(1-\zeta_5), 2^2 \pi_2, 2^4 \pi_3} \cap \lara{2^4\pi_1, 2^2(1-\zeta_5), 2\pi_2, 2^2\pi_3}.
% \end{align*}
% Arbitrary element of each space can be written by
% \begin{align*}
%     2^{a + b+ 2c + 4d} \pi_1^a (1-\zeta_5)^b \pi_2^c \pi_3^d, \qquad
%     2^{4a + 2b + c + 2d}\pi_1^a (1-\zeta_5)^b \pi_2^c \pi_3^d.
% \end{align*}
% Hence finding the solution space of $a+b+2c+4d = 4a+2b+c+2d$, which is generated by $(1,0,1,1), (0,1,1,0), (0,2,0,1)$, we get three more vectors $2^2 \pi_1  \pi_2 \pi_3, 2^3(1-\zeta_5)\pi_2, 2(1-\zeta_5)^2 \pi_3$.
% Now, we can do the same thing for $V_3 \cap V_4$, and we have
% \begin{align*}
%     V_3 \cap V_4 = \lara{ \zeta_5, 2^4 \pi_2, 2^3 \pi_4(1-\zeta_5), \pi_3(1 + \zeta_5)^4, 2\pi_1(1-\zeta_5), 2^2 \pi_1\pi_3\pi_4}.
% \end{align*}
% Note that $2(1+\zeta_5) \pi_3(1+\zeta_5)^4 \in V_1\cap V_2 + V_3 \cap V_4$, so we have 
% $V_1 \cap V_2 + V_3 \cap V_4$ is $K(S, 5)$.
% Hence, we know that the dimension of $V_1 \cap V_2 \cap V_3 \cap V_4$ is 4.
One can check that 
\begin{align*}
    W:=V_1 \cap V_2 \cap V_3 \cap V_4 = \lara{ \zeta_5, 2\pi_1\pi_2\pi_3\pi_4, 2^2\pi_2\pi_4(1-\zeta_5), 2^4(1-\zeta_5)^2(1+\zeta_5)^4\pi_1\pi_3\pi_4^3}.
\end{align*}

We recall that our embedding of $K$ into $K_5$ maps $\zeta_5$ to $1 + \pi_5$.
Then, $-\pi_1, \pi_2, \pi_3, -\pi_4$ are also maps to \begin{align*}
    \begin{array}{ll}
    \pi_1 &\mapsto -(1 + 3\pi_5 + 4\pi_5^2 + \pi_5^3 + \pi_5^4)  \\
    \pi_2 &\mapsto 1 + \pi_5 + 3\pi_5^2 + 2\pi_5^3 + 3\pi_5^4 + 4\pi_5^5  \\
    \pi_3 &\mapsto 1 + 2\pi_5 + \pi_5^2 + 4\pi_5^3 + 2\pi_5^4 \\
    \pi_4 &\mapsto -(1 + 4\pi_5 + 2\pi_5^2 + 3\pi_5^3 + \pi_5^5) 
    \end{array}
\end{align*}
modulo $O(\pi_5^6)$, which corresponds to the $U^{(1)}$-part. By a routine computation, we have
\begin{align*}
    i_5(\pi_1, \pi_2, \pi_3, \pi_4) = (
    \beta^3 \gamma \delta^2 \epsilon^2 \eta^3,
    \beta \gamma^3 \delta^4 \epsilon \eta^3,
    \beta^2 \delta^4 \epsilon^4 \eta^2,
    \beta^4 \gamma \epsilon^3 \eta^2).
\end{align*}
We already know that 
\begin{align*}
    i_5(2, \zeta_5, 1 + \zeta_5, 1- \zeta_5)
    = (\epsilon^3 \eta, \beta\gamma\epsilon, \beta^2\gamma^4\delta^2\epsilon^4, \alpha)
\end{align*}
and $F_5(J_{m^2}(K_5)/ \phi J_{m^2}(K_5)) = \lara{\delta, \epsilon, \eta}$ by Proposition \ref{prop:map F on 5}.
The images of our basis members of $W$ in the quotient space $K_5^{\times}/K_5^{\times 5}\cdot F_5(J_{m^2}(K_5)/ \phi J_{m^2}(K_5))$ are $\overline{\beta\gamma}, \overline{1}, \overline{\alpha \gamma^4}, \overline{\alpha^2}$, respectively. 
Therefore $\Sel_{\phi}(J_{101^2}/\bQ)$ is one dimensional vector space generated by $2\pi_1\pi_2\pi_3\pi_4$.
\end{proof}

We conclude this section with an example on general $m$ which is not divided by a rational prime equivalent to one modulo five.

\begin{example}[$m = p_1p_2$ where $(p_1, p_2) \equiv (3, 4) \pmod{25}$.]
\label{exam:}
Let $p_1 \equiv 3$ and $p_2 \equiv 4 \pmod{25}$, and $\pi_w$ and $\pi_{\overline{w}}$ be prime elements $a \pm b \sqrt{5}$ for $a, b\in \frac{1}{2}\bZ$ of $K$ lying over $p_2$.
Then, by Remark \ref{rem:} and Lemma \ref{lem:localover5},
\begin{equation*}
    F_5(J(K_5) / \phi(J(K_5))) = \lara{\gamma \delta, \epsilon, \eta} \textrm{ and }
    \im i_5 = \lara{\beta\gamma\epsilon, \beta^2\gamma^4\delta^2\epsilon^4, \epsilon^3\eta, \alpha, \epsilon\eta^2, (\gamma\delta\epsilon)^b}.
\end{equation*}
So the previous argument shows that
\begin{equation*}
    (i_5^{-1} \circ F_5) (J(K_5)/\phi J(K_5))
    = \left\{ \begin{array}{cc}
    \lara{2, p_1}     &  \textrm{if } b \not\equiv 0 \pmod{5}, \\
    \lara{2, p_1, \pi_w, \pi_{\overline{w}}}     &  \textrm{if } b\equiv 0 \pmod{5}.
    \end{array} \right.
\end{equation*}
For the other bad primes $\fp$ we have $(i_\fp^{-1}\circ F_\fp)(J(K_\fp) / \phi J(K_\fp))$ contains $\lara{2, p_1, \pi_w, \pi_{\overline{w}}}$, by Proposition \ref{prop:descent not over 5}.
Therefore, for such $m = p_1p_2$,
\begin{equation*}
    \dim_{\bF_5}\Sel_\phi(J_m/\bQ) = \left\{
    \begin{array}{cc}
    2     & \textrm{if } b \not\equiv 0 \pmod{5}, \\
    4     & \textrm{if } b \equiv 0 \pmod{5}.
    \end{array}
    \right.
\end{equation*}
\end{example}

% \textcolor{blue}{Totally split?}
% Also, the totally split case can be dealt with the similar way. The main difficulty of totally split cases is that there are too many cases to describe the image of the uniformizer in $K_5^\times/K_5^{\times 5}$.
% However, for a fixed prime $p$, it is not so hard to do a descent for $J_{p^2}$ even when $p$ is totally split in $K/\bQ$.
% \begin{example}[$m = 191$]
% (or simpler prime)
% \end{example}

\section{Special values of $L$-functions} \label{sec:L-value}

In this section we will find sufficient conditions on $A$ such that $L(1,J_A)$ becomes nonzero.
By \cite[Theorem 4]{Mil72}, there is a Hecke character $\eta_A$ of $K$ such that
\begin{align*}
L(s,J_A)=L(s,\eta_A).
\end{align*}
Following \cite[section 2]{SY03}, we denote $F:=\mathbb{Q}(\sqrt{5})$ and $\chi_A:=\eta_A|\cdot|_\mathbb{A}^{1/2}$ with $\mathbb{A}:=\mathbb{A}_F$ the ring of ad\`eles so that
\begin{align*}
L(1,J_A)=L(1,\eta_A)=L\left(\frac{1}{2},\chi_A\right).
\end{align*}
From now on, we assume that the global root number of $\chi_A$ is $1$. 
Based on the work of \cite{Yan97, Yan99}, Stoll and Yang give the following:
\begin{proposition}[{\cite[Proposition 3.1]{SY03}}]
\label{prop:Stoll-Yang3.3}
With the notation in \cite{SY03}, we have
\begin{align*}
    L(1, \eta_A) = \frac{\pi^2}{50C_1C_2} 
    \left| \sum_{x \in F} \prod_{v \nmid 2A} \phi_v(x) \prod_{v \mid 2A} I_v(x) \right|^2
\end{align*}
for some constant $C_1$ and $C_2$.
\end{proposition}
Here $\phi=\prod_v\phi_v\in S(\mathbb{A})$ is an appropriately chosen Schwartz-Bruhat function and 
\begin{align} \label{eqn:defIv}
	I_v(x) = \int_{G_v} \omega_{\alpha, \chi_A, v }(g) \phi_v(x) dg
\end{align}
as in \cite[p. 277]{SY03}.
%$I_v$ is a function given by integrating the action of a Weil representation over some domain.
We will introduce more precise notations later.
Stoll and Yang further give a concrete choice of $\phi_v$ for $v \nmid A$ and infinite $v$. It allows them to compute $L(1, \eta_1)$.
In this paper, we choose $\phi_v$ for $v \mid A$ and consider when $I_v(x)$ is non-zero.

Since the global root number of $\chi_A$ is $+1$, there is a unique $\alpha\in F^\times$ up to norm from $K^\times$ such that 
\begin{align*}
\prod_{\substack{\textrm{$w$ places of $K$} \\ w\mid v}}\epsilon\left(\frac{1}{2},\chi_{A,w},\frac{1}{2}\psi_{K_w}\right)\chi_{A,w}(\delta)=\epsilon_v(\alpha)
\end{align*}
for all places $v$ of $F$ (cf. \cite[p. 276]{SY03}).
Here $\delta:=\zeta_5^{-2}-\zeta_5^2$, 
$\psi$ is an additive character of $\mathbb{A}_F$ given by $\psi = \prod_v \psi_v$ for $\psi_v(x) = e^{-2\pi \sqrt{-1}\lambda_v(x)}$ where 
\begin{align*}
\xymatrix{\lambda_v : F_v \ar[r]^-{\mathrm{Tr}_{F_v/\mathbb{Q}_p}} & \bQ_p \ar[r] & \bQ_p/\bZ_p \ar[r] & \bQ/\bZ},
\end{align*}
and $\psi_K := \psi \circ \mathrm{Tr}_{K/F}$.
Also, $\epsilon$ on the left hand side are the local root numbers as in \cite[Proposition 2.2]{SY03}, and $\epsilon_v$ is the local part of the Hecke character belonging to $K/F$. We let rings act on additive characters defined on them by multiplication with arguments. For example,
\begin{align*}
\left(\frac{1}{2}\psi_{K_w}\right)(x):=\psi_{K_w}\left(\frac{1}{2}x\right).
\end{align*}
Since we only concern the case where $A$ is a square not divisible by $2$, \cite[Lemma 2.3]{SY03} tells us that we may choose
\begin{align*}
\alpha\in\left(\prod_{2\neq p|A}p\right)\cdot N_{K/F}K^\times
\end{align*}
where $N_{K/F}$ denotes the norm. Next, we need to choose an appropriate Schwartz-Bruhat function $\phi=\prod_v\phi_v\in S(\mathbb{A})$ as in \cite[p. 277]{SY03}.
To be more precise, we introduce more notations in \cite[section 2]{SY03}. We fix an embedding $K\hookrightarrow\mathbb{C}$ such that $\zeta_5\mapsto\exp(2\pi\sqrt{-1}/5)$. 
We also fix a CM type $\Phi=\{\sigma_2,\sigma_4\}$ of $K$ where $\sigma_r(\zeta_5)=\exp(2\pi r\sqrt{-1}/5)$.  Then the following lemma tells us a possible choice of $\phi_v$ for almost all places $v$.
\begin{lemma}[{\cite[Lemma 3.2]{SY03}}] \label{lem:Stoll-Yang 3.2}
Denote $\mathrm{char}(X)$ the characteristic function of the set $X$.
Then,
\begin{align*}
\phi_v(x)=\left\{\begin{array}{ll}
\mathrm{char}(\mathcal{O}_{F,v})  & v\nmid10A\infty,\alpha\in\mathcal{O}_{F,v}^\times , \\
|2\sigma_j(\alpha\delta^3)|^{1/4}\exp\left(-\pi|\sigma_j(\alpha\delta^3)|\sigma_j(x)^2\right) & v=\sigma_j\in\{\sigma_2,\sigma_4\} .
\end{array}\right.
\end{align*}
\end{lemma}
If we choose $\alpha\in F^\times$ as above such that $\alpha\in\mathbb{Z}_2^\times$, then \cite[Corollary 5.8]{SY03} tells us that we may choose
\begin{align*}
\phi_2=\mathrm{char}\left(\frac{1}{2}+\mathcal{O}_{F,2}\right).
\end{align*}
We note that $\phi_2 = I_2$ and $I_2$ is a constant function (See \cite[\S 4]{SY03}). 
%\textcolor{blue}{and we should note that $\phi_2 = I_2$, and $I_2$ is a constant function? (written in \cite[\S 4]{SY03}) }
At $v=\sqrt{5}$, \cite[Proposition 1.2, Corollary 1.4]{Yan99} tell us that we may choose
\begin{align*}
\phi_{\sqrt{5}}=5^{\frac{2n(\chi_{A,\lambda})-1}{4}}\xi_{\lambda}\cdot\mathrm{char}(\mathcal{O}_{F,\sqrt{5}}).
\end{align*}
Here, by denoting $\Delta := \delta^2$,
\begin{quote}
(1) $\lambda:=1-\zeta_5\in K$ is a prime element lying over $\sqrt{5}$.\\
(2) $n(\chi_{A,\lambda})$ is the conductor exponent of $\chi_{A,\lambda}$ which is completely determined by $q_5(A)=(A^4-1)/5$ (see \cite[Proposition 2.2 (5)]{SY03}):
\begin{align*}
n(\chi_{A,\lambda})=\left\{\begin{array}{ll}
1 & \textrm{if}\quad5\mid q_5(A), \\
2 & \textrm{if}\quad5\nmid q_5(A).
\end{array}\right.
\end{align*}
(3) With $G=\{\pm1\}\times U_K^{(1)}$, write $g=x+y\delta\in G$ and set
\begin{align*}
\xi_{\lambda}(g)=\left\{\begin{array}{ll}
\chi_{A,\lambda}(\delta(g-1))(\Delta,-y)_F & \textrm{if}\quad g\in U_K^{(1)}, \\
\chi_{A,\lambda}(\delta(g-1))(\Delta,-2\alpha)_F\epsilon(\frac{1}{2},\epsilon_{K_w/F_v},\psi_{K_\lambda}) & \textrm{if}\quad g\in G\setminus U_K^{(1)}.
\end{array}\right.
\end{align*}
This comes from \cite[Proposition 1.2 (1)]{Yan99}.\footnote{
It seems that there is a typo in \cite[Proposition 1.2 (1)]{Yan99}.
Compare the statement and its proof \cite[pp.354-455]{Yan99}.}
\end{quote}
By Proposition \ref{prop:Stoll-Yang3.3} and Lemma \ref{lem:Stoll-Yang 3.2}, we obtain
\begin{align} \label{eqn:L prever}
	L(1, \eta_A) =  \frac{\pi^2}{50C_1C_2} \cdot  5^{\frac{2n(\chi_{A,\lambda})-1}{4}} \cdot \left|
	\sum_{x \in X_A'}\xi_{\lambda}(x) \phi_{\sigma_2}(x) \phi_{\sigma_4}(x) \cdot \lbrb{ \prod_{v \mid A} I_v(x) }
	\right|
\end{align}
where 
\begin{align*}
	X_A' = F \cap \lbrb{ \bigcap_{ v\nmid 2A\infty  }\cO_{F, v} } \cap \lbrb{ \frac{1}{2} + \cO_{F, 2} }.
\end{align*}

For $v\mid A$ and $w$ a place of $K$ dividing $v$, we always have $n(\chi_{A,w})=1$ by \cite[Proposition 3.3]{Sto2}. 
First, we consider the case $v\mid A$ splits in $K/F$. In this case we apply \cite[section 2]{Yan97}. Under the identification
\begin{align*}
K_v\cong\frac{F[t]}{(t^2-\Delta)}\otimes_FF_v\cong F_v\cdot\delta\oplus F_v\cdot(-\delta)
\end{align*}
we have $\delta=(1,-1)\in F_v\oplus F_v$. 
Denote $\pi_{F_v}\in\mathcal{O}_{F,v}$ by a uniformizer and in this case $n_v = 1.$ To get $\phi_v=\phi_{v,1}$, following the notation of \cite[Theorem 2.15]{Yan97}, we first compute
\begin{align*}
\rho\left(\mathrm{char}\left(1+\pi_{F_v} \mathcal{O}_{F,v}\right)\right)(x)&:=
|\alpha|_v^{\frac{1}{2}}\psi_v\left(\frac{\alpha x^2}{2}\right)\int_{F_v}\psi_v(\alpha xy)\psi_v\left(\frac{\alpha y^2}{4}\right)\mathrm{char}\left(1+\pi_{F_v}\mathcal{O}_{F,v}\right)(y)dy\\
&=|\alpha|_v^{\frac{1}{2}}\psi_v\left(\frac{\alpha x^2}{2}\right)\int_{1+\pi_{F_v}\mathcal{O}_{F,v}}\psi_v(\alpha xy)dy\\
&=|\alpha|_v^{\frac{1}{2}}\psi_v\left(\frac{\alpha x^2}{2}\right)\int_{\pi_{F_v}\mathcal{O}_{F,v}}\psi_v(\alpha x(y+1))dy\\
&=|\alpha|_v^{\frac{1}{2}}\psi_v\left(\frac{\alpha x^2}{2}+\alpha x\right)\int_{\pi_{F_v}\mathcal{O}_{F,v}}\psi_v(\alpha xy)dy\\
&=|\alpha|_v^{\frac{1}{2}}\psi_v\left(\frac{\alpha x^2}{2}+\alpha x\right)\mathrm{meas}(\pi_{F_v}\cO_{F, v}) \mathrm{char}\left(\pi_{F_v}^{-2}\mathcal{O}_{F,v}\right)(x).
\end{align*}
Hence we get
\begin{align*}
\phi_v=\mathrm{meas}(\mathcal{O}_{F,v})^{-\frac{1}{2}}\mathrm{meas}(\pi_{F_v}\cO_{F, v}) q_v^{\frac{1}{2}}|\alpha|_v^{\frac{1}{2}}\psi_v\left(\frac{\alpha x^2}{2}+\alpha x\right)\mathrm{char}\left(\pi_{F_v}^{-2}\mathcal{O}_{F,v}\right)(x).
\end{align*}
To apply \cite[Proposition 3.1]{SY03}, we need to compute
\begin{align*}
I_v(x)&:=\int_{\mathcal{O}_{F,v}^\times}\omega_{\alpha,\chi_A,v}(g)\phi_v(x)dg\\
&=\int_{\mathcal{O}_{F,v}^\times}\chi_{A,v}(g)|g|_v^{\frac{1}{2}}\phi_v(xg)dg\\
&=\int_{\mathcal{O}_{F,v}^\times}\phi_v(xg)dg\\
&=\mathrm{meas}(\mathcal{O}_{F,v})^{-\frac{1}{2}}\mathrm{meas}(\pi_{F_v}\cO_{F, v}) q_v^{\frac{1}{2}}|\alpha|_v^{\frac{1}{2}}
\int_{\mathcal{O}_{F,v}^\times}\psi_v\left(\frac{\alpha}{2}(xg)^2+\alpha(xg)\right)\mathrm{char}\left(\pi_{F_v}^{-2}\mathcal{O}_{F,v}\right)(xg)dg\\
&=\mathrm{meas}(\mathcal{O}_{F,v})^{-\frac{1}{2}}\mathrm{meas}(\pi_{F_v}\cO_{F, v}) q_v^{\frac{1}{2}}|\alpha|_v^{\frac{1}{2}}\mathrm{char}\left(\pi_{F_v}^{-2}\mathcal{O}_{F,v}\right)(x)\int_{\mathcal{O}_{F,v}^\times}\psi_v\left(\frac{\alpha}{2}(xg)^2+\alpha(xg)\right)dg.
\end{align*}
We note that the action of Weil representation $\omega$ is described in \cite[Corollary 2.10]{Yan97}.
Since there is a representative 
\begin{align*}
\alpha\in\left(\prod_{2\neq p|A}p\right)\cdot N_{K/F}K^\times,
\end{align*}
we choose $\alpha$ such that $\psi_v\left(\frac{\alpha}{2}(xg)^2+\alpha(xg)\right) = 1$ for $g \in \cO_{F, v}^\times$ and $x \in \pi_{F_v}^{-2} \cO_{F, v}$ for all $v \mid A$ splitting in $K/F$.
%$\ord_{F_v}(\alpha) \geq n_v - 2$ for all $v \mid A$ splitting in $K/F$.
Then
\begin{align}\label{valueIv}
I_v|_{\pi_{F_v}^{-2}\mathcal{O}_{F,v}}=\mathrm{meas}(\mathcal{O}_{F,v})^{-\frac{1}{2}} \mathrm{meas}(\pi_{F_v}\cO_{F, v}) q_v^{\frac{1}{2}}|\alpha|_v^{\frac{1}{2}}\int_{\mathcal{O}_{F,v}^\times}dg
=\frac{\mathrm{meas}(\mathcal{O}_{F,v}^\times)}{\mathrm{meas}(\mathcal{O}_{F,v})^{\frac{1}{2}}}
\mathrm{meas}(\pi_{F_v}\cO_{F, v}) q_v^{\frac{1}{2}}|\alpha|_v^{\frac{1}{2}}
\end{align}
is a non-zero constant. Therefore, there is a non-zero constant $c_v(\alpha)$ such that 
\begin{align} \label{eqn:Ivsplit}
	I_v(x) = c_v(\alpha) \mathrm{char} (\pi_{F_v}^{-2}\cO_{F, v}),
\end{align}
when $v \mid A$ splits in $K/F$.

Finally, consider the case $v\mid A$ is inert in $K$. Following the notation of \cite[p. 339]{Yan99}, we have
\begin{align*}
n(\psi_{K_v}')=n\left(\frac{\alpha\delta}{4}\psi_{K_v}\right)=n(\psi_{K_v})-\mathrm{ord}_{F_v}(\alpha)=-\mathrm{ord}_{F_v}(\alpha).
\end{align*}
We choose $\alpha$ so that $\mathrm{ord}_{F_v}(\alpha)=1$ and $n(\psi_{K_v}')=-1$.
Since we have $n(\chi_{A,v})=1$ and $w\mid v$ is unramified, we are in the case of \cite[Proposition 1.5]{Yan99} with $\eta=1$ the trivial character. 
Then we may choose, 
\begin{align} 
\phi_v(x)&=\mathrm{char}(\pi_{F_v}\mathcal{O}_{F,v})(\pi_{F_v}x)+ \nonumber \\
&\quad+\frac{1}{2G(\psi_{F_v}'')}\sum_{\substack{(S,T)\in\kappa_v^2 \\ S^2-T^2\equiv\Delta\bmod\pi_{F_v}}}\xi_v^{-1}\left(\frac{S+\delta}{T}\right)\left(\frac{T}{\kappa_v}\right)\psi_{F_v}''\left(\frac{\Delta\alpha}{2}S(\pi_{F_v}x)^2\right)\mathrm{char}(\mathcal{O}_{F,v})(\pi_{F_v}x) \label{eqn:phiv1}
\end{align}
when $\xi_v(-1) = \quadsym{-1}{\kappa_v}$,
or %\footnote{My choice for the second case is $a=1$ in \cite[Proposition 1.5 (2)]{Yan99} one may use other choices to simplify computation if possible}
\begin{align}
\phi_v(x)&:=\mathrm{char}(1+\pi_{F_v}\mathcal{O}_{F,v})(\pi_{F_v}x)-\mathrm{char}(-1+\pi_{F_v}\mathcal{O}_{F,v})(\pi_{F_v}x)+ \nonumber \\
&\quad+\frac{1}{G(\psi_{F_v}'')}\sum_{\substack{(S,T)\in\kappa_v^2 \\ S^2-T^2\equiv\Delta\bmod\pi_{F_v}}}\xi_v^{-1}\left(\frac{S+\delta}{T}\right)\left(\frac{T}{\kappa_v}\right)\psi_{F_v}''(S(\pi_{F_v}x)^2-2T\pi_{F_v}x+S)\mathrm{char}(\mathcal{O}_{F,v})(\pi_{F_v}x)  \label{eqn:phiv2}
\end{align}
when $\xi_v(-1) = - \quadsym{-1}{\kappa_v}$ and $\xi_v^{-1} \neq \eta_0$,
where $\kappa_v:=\mathcal{O}_{F,v}/\pi_{F_v}$ is the residue field of $F_v$. Note that $\psi_{F_v}''$ in \cite[Proposition 1.5]{Yan99} has conductor $\pi_{F_v}\mathcal{O}_{F,v}$ (see the proof of \cite[Proposition 3.4]{Yan98} for the detail) so we regard $\psi_{F_v}''$ as a character of $\kappa_v$ and $G(\psi_{F_v}'')$ is the Gauss sum of $\psi_{F_v}''$.
%We define
%\begin{align} \label{eqn:defIv}
%	I_v(x) = \int_{G_v} \omega_{\alpha, \chi_A, v }(g) \phi_v(x) dg
%\end{align}
%as in \cite[p.277]{SY03}.
Together with (\ref{eqn:L prever}), we obtain
\begin{proposition} \label{prop:Lvalue}
Suppose that the root number of $\eta_A$ is $+1$. Then, there is a non-zero constant $c_v(\alpha)$ such that 
\begin{align} \label{eqn:Lvalue equation}
	L(1, \eta_A) = \frac{\pi^2}{50C_1C_2}  \cdot  5^{\frac{2n(\chi_{A,\lambda})-1}{4}} \cdot  \prod_{\substack{v \mid A \\ v \textrm{ split} }} c_v(\alpha) 
	\cdot\left|
	\sum_{x \in X_A} \xi_{\lambda}(x) \phi_{\sigma_2}(x) \phi_{\sigma_4}(x) \cdot 
	 \prod_{\substack{v \mid A \\ v \textrm{ inert} }} I_v(x) \right|
\end{align}
where $I_v(x)$ is taken from (\ref{eqn:defIv}) and
\begin{align*}
	X_A = F \cap \lbrb{ \bigcap_{ \substack{v\nmid 2A\infty } }\cO_{F, v} } \cap \lbrb{ \frac{1}{2} + \cO_{F, 2} } \cap \lbrb{\bigcap_{\substack{v \mid A \\ v \textrm{ split} }} \pi_{F_v}^{-2}\cO_{F, v}}.
\end{align*}
\end{proposition}

\begin{proof}[Proof of Theorem \ref{thm:main2}]
When $5^2 \mid (A^4-1)$, we have $n(\chi_{A, \lambda}) = 1$ which implies that $\xi_\lambda$ is trivial (See \cite[Proposition 1.2, Corollary 1.4]{Yan99}).
Since every prime divisor of $A$ splits in $K/F$, we obtain that
\begin{align} \label{eqn:Lvalueproof}
	L(1, \eta_A) = \frac{\pi^2}{50C_1C_2} \cdot 5^{\frac{1}{4}} \cdot \prod_{\substack{v \mid A \\ v \textrm{ split} }} c_v(\alpha) 
	\cdot\left|
	\sum_{x \in X_A} \phi_{\sigma_2}(x) \phi_{\sigma_4}(x) \right|.
\end{align}
Recall that $\phi_2$ and $\phi_4$ have real values on $F$. Therefore,
\begin{align*}
\phi_{\sigma_2}(x)\phi_{\sigma_4}(x)=\sqrt{2}\alpha^{\frac{1}{2}}5^{\frac{3}{8}}\exp\left(-\pi\alpha\left(\left(2\sin\frac{2\pi}{5}\right)^3\sigma_2(x)^2+\left(2\sin\frac{4\pi}{5}\right)^3\sigma_4(x)^2\right)\right)
\end{align*}
is non-zero and the last term of (\ref{eqn:Lvalueproof}) does not vanish. Hence $L(1, \eta_A)$ is non-zero.
\end{proof}

\begin{proof}[Proof of Corollary \ref{cor:101}]
We note that $q_5(101^2)$ is divided by 5.
Now the result follows from Propositions \ref{prop:101} and Theorem \ref{thm:main2}. 
\end{proof}

\emph{Acknowledgement}
Authors thank to Dohyeong Kim for the useful suggestions.
K. Jeong is supported by the National Research Foundation of Korea (NRF) grant funded by the Korea government (MSIT) 
(2019R1C1C1004264 and 2020R1A4A1016649).
D. Yhee is supported by the National Research Foundation of Korea(NRF) Grant funded by the Korean Government (MSIT) (2017R1A5A1015626).

\end{document}